\documentclass[letterpaper, 10 pt, conference]{ieeeconf}
\IEEEoverridecommandlockouts  
\usepackage{amsmath,amssymb,amsfonts,mathrsfs,soul,bbm,tikz}
\usetikzlibrary {shadings,fadings,arrows.meta}
\usepackage{booktabs} 
\usepackage{siunitx} 
\usepackage{hyperref}

\usepackage{amsthm} 
\usepackage{sansmath}
\newtheorem{thm}{Theorem}

\newtheorem{assum}{Assumption}
\newtheorem{prop}{Proposition}
\newtheorem{defn}{Definition}
\newtheorem{problem}{Problem}

\usepackage{graphicx}
\usepackage{subcaption}
\graphicspath{{../Figures/}}
\usepackage{listings}
\usepackage{mathtools}
\usepackage{algorithm}
\usepackage{algpseudocode}
\usepackage{multirow}
\usepackage[labelsep=period]{caption}
\usepackage{upgreek}
\usepackage{cite}
\usepackage{lipsum}

\newcommand{\ind}[1]{{1}#1}

\renewcommand{\d}{\mathrm{d}}

\algrenewcommand\algorithmicrequire{\textbf{Input:}}
\algrenewcommand\algorithmicensure{\textbf{Output:}}

\usepackage{xcolor}

\title{Stochastic Reachability of Uncontrolled Systems via Probability Measures:\\
Approximation via Deep Neural Networks}
\author{Karthik Sivaramakrishnan, Vignesh Sivaramakrishnan, Rosalyn A. Devonport, Meeko M. K. Oishi
    \thanks{
    This material is based upon work supported by the National Science Foundation under NSF Grant Number CNS-1836900.  Any opinions, findings, and conclusions or recommendations expressed in this material are those of the authors and do not necessarily reflect the views of the National Science Foundation.  The NASA University Leadership initiative (Grant \#80NSSC20M0163) provided funds to assist the authors with their research, but this article solely reflects the opinions and conclusions of its authors and not any NASA entity.  This material is based upon research supported, in collaboration with Verus Research, by the Air Force Research Lab (AFRL) under agreement number FA9453-23-C-A025. The U.S. Government is authorized to reproduce and distribute reprints for Governmental purposes notwithstanding any copyright notation thereon.
    The views and conclusions contained herein are those of the authors and should not be interpreted as necessarily representing the official policies or endorsements, either expressed or implied, of Air Force Research Laboratory (AFRL) and or the U.S. Government.
    }
    \thanks{Karthik Sivaramakrishnan is a graduate student with Electrical and Computer Engineering, University of New Mexico, Albuquerque, NM.\newline Email: \tt{ksivaramakrishnan20@unm.edu}.}
    \thanks{Vignesh Sivaramakrishnan is a graduate student with Electrical and Computer Engineering, University of New Mexico, Albuquerque, NM.\newline Email: \tt{vigsiv@unm.edu}.}
    \thanks{Rosalyn Alex Devonport is a postdoctoral scholar with Electrical and Computer Engineering, University of New Mexico, Albuquerque, NM.\newline Email: \tt{devonport@unm.edu}.}
    \thanks{Meeko M. K. Oishi is a professor with Electrical and Computer Engineering, University of New Mexico, Albuquerque, NM.\newline Email: \tt{oishi@unm.edu}. }
}
        
\begin{document}
\maketitle
\begin{abstract}
    This paper poses a theoretical characterization of the stochastic reachability problem in terms of probability measures, capturing the probability measure of the state of the system that satisfies the reachability specification for all probabilities over a finite horizon. 
    We achieve this by constructing the level sets of the probability measure for all probability values and, since our approach is only for autonomous systems, we can determine the level sets via forward simulations of the system from a point in the state space at some time step in the finite horizon to estimate the reach probability.
    We devise a training procedure which exploits this forward simulation and employ it to design a deep neural network (DNN) to predict the reach probability provided the current state and time step. 
    We validate the effectiveness of our approach through three examples. 
\end{abstract}

\section{Introduction}

Stochastic reachability is an important tool to provide probabilistic assurances of safety in stochastic, dynamical systems, by ensuring that constraints on the state space are met with at least a desired likelihood, despite bounded control authority.  However, computing stochastic reachable sets is a difficult challenge, because the solution to the stochastic reachability problem is based in dynamic programming \cite{abate2008probabilistic,summers2010verification,SUMMERS20132906}.  While some progress has been made by exploiting linearity in the dynamics and structure in the stochasticity \cite{lesser2013stochastic,vinod2017scalable,vinod2021stochastic, kariotoglou2017linear}, computational solutions remain challenging, particularly for systems described by general nonlinear dynamics or with arbitrary stochastic processes.  In such cases, methods based in approximate dynamic programming have been proposed \cite{kariotoglou2013approximate,thorpe2022state}, that learn the value function at each time step, which can be expensive for longer horizons.

\textit{The main contribution of this paper is a formulation of the stochastic reachability problem via probability measures, that captures the state distribution which satisfies the reachability specification at every time step.}
We synthesize a probability measure that ensures compliance with reachability specifications for all probabilities.   
We propose a measure-based approach to approximately computing reachable sets using neural nets as surrogate functions. 
The probability measure of a random variable returns, for a set input (e.g., a set of states), the likelihood of the variable lying in that set.  With a probability measure, we can find a set associated with any prescribed likelihood~\cite{bertsekas1996stochastic, chow1997probability}.
Occupation measures, a type of probability measure, have been used to compute regions of attraction \cite{henrion2013convex} and control invariant sets~\cite{korda2014convex}. Essentially, we use measures to represent the stochastic nature of the evolution of the system state. 
We additionally seek numerical implementations that can approximate stochastic reachable sets by exploiting this formalism.  We capture the level sets for some likelihood as a root-finding problem over possible probability measures, then propose a data-driven approach based in empirical assessments of constraint satisfaction. This approach reduces problem of constructing a state probability measure to one of forward propagation of the dynamics, to empirically evaluate probabilities. 
We construct these measures using a backwards recursion formalism; this formalism yields a sequence of value functions that we estimate from forwards trajectory data using a neural network surrogate function. The value functions themselves are not novel-- indeed, they are central to the dynamic programming approach of~\cite{abate2008probabilistic}-- but they are intensive to compute in general. 
The method by which we approximate them from data, and to use them to generalize the observed evolutions into an approximate reachable set over the continuum of the state space, is another of our contributions. This approach to approximate reachability has proven successful in other domains where dynamic programming is employed, e.g. in~\cite{jiang2017usingneuralnetworkscompute}, but has not yet been applied to the case of stochastic reachability. We focus on {\em autonomous systems}, that is, dynamical systems with no control input.  

Our deep neural net learns the probability measure across all safety likelihoods simultaneously, but has computational complexity that scales with the network size (as determined by the number of neurons and layers).

In the following, 
Section~\ref{sec:prelim} presents the preliminaries and the problem formulation.  In Section \ref{sec:theory}, 
We characterize the stochastic reachability problem via probability measures. 
Section \ref{sec:NN}
describes the training procedure to create a deep neural net that predicts the stochastic reachability probability for a finite time horizon, for a given state and time step. 
Finally, in Section \ref{sec:example}, we apply our approach to multiple scenarios, including a satellite pointing problem.

\section{Preliminaries and Problem Formulation}
\label{sec:prelim}
We presume Borel space, $(\mathcal{X},\mathscr{B}(\mathcal{X}))$ where $\mathcal{X}\subseteq\mathbb{R}^n$ and $\mathscr{B}(\mathcal{X})$ is a Borel sigma algebra, i.e. a non-empty collection of subsets $X\in\mathscr{B}(\mathcal{X})$~\cite{bertsekas1996stochastic}. 
A probability measure, $\mu: \mathscr{B}(\mathcal{X})\rightarrow[0,1]$, maps from the Borel sigma algebra outputs a probability value between zero and one. 
Random variables are in boldface, $\mathbf{x}$, where we specify the probability measure we define them over as $\mathbf{x}\sim\mu$.

\begin{defn}[{\cite[Definition 7.12]{bertsekas1996stochastic}}]
    \label{def:TransitionKernel}
    A discrete-time, autonomous stochastic system is a Borel measurable transition kernel, $Q: \mathscr{B}(\mathcal{X}) \times \mathcal{X}\rightarrow [0,1]$, that takes a set from the Borel sigma algebra and outputs a value between zero and one. 
    We denote this function as
    \begin{equation}
        \label{eq:transitionKernel}
        Q(\cdot|x).
    \end{equation}
\end{defn}
We denote whether a state resides within a set, $T\in\mathscr{B}(\mathcal{X})$, by an indicator function,
\begin{equation}
    \label{eq:indicator}
    \ind_{T}(x) = 
    \begin{cases}
        1,\ x\in T\\
        0,\ x\notin T. 
    \end{cases}
\end{equation}
The probability of the current state residing within a set, provided we know the probability measure of the current state, $\mu_k(\cdot)$, is described by  the integral
\begin{equation}
    \label{eq:StateInSet}
    \mu_{k}(T_{k}) = \int_{\mathcal{X}}\ind_{T_{k}}(x_{k})\mu_k(\mathrm{d}x_k). 
\end{equation}
In addition, we can relate the probability measure of the next state for all sets in the Borel sigma algebra to the probability measure of the current state by
\begin{equation}
    \label{eq:measureState}
    \mu_{k+1}(X) = \int_{\mathcal{X}}Q(X|x_k)\mu_k(\mathrm{d}x_k), \: \forall X\in\mathscr{B}(\mathcal{X})
\end{equation}
We note that $\mu_{k+1}$ is unique \cite[Lemma 10.4.3]{bogachev2007measure}.

\begin{problem}
    \label{prob:RAviaPM}
    Given a target tube, $T_k$ for $k=\{0,\hdots,N\}$, find the state probability measures $\mu_k(\cdot)$ for a finite horizon $N$ such that the system in Definition~\ref{def:TransitionKernel} resides in the target tube for all probability values $\alpha\in[0,1]$.  
\end{problem}
Problem~\ref{prob:RAviaPM} can be adapted to different types of stochastic reachability specifications including reach-avoid~\cite{kariotoglou2013approximate}, terminal-hitting time~\cite{abate2008probabilistic}, and, as it is now, for target-tube~\cite{vinod2021stochastic}.
This is accomplished by writing the reachability specifications using the appropriate indicator functions in ~\eqref{eq:indicator}. 

We address Problem~\ref{prob:RAviaPM} by providing a theoretical characterization of the solution to Problem~\ref{prob:RAviaPM} in Section~\ref{sec:theory}, and a numerical implementation based in deep neural nets (DNNs) in Section~\ref{sec:NN}.

\section{The Measure-Theoretic Characterization of Stochastic Reachability}
\label{sec:theory} 

We now present the measure-theoretic formulation of the stochastic reachability problem that will form the basis for our neural network-based reachable set approximation method.
First, we define the stochastic reachability specifications in terms of probability measures; then, we characterize how to use measures to solve stochastic reachability problems.

\subsection{Backward Recursion of the State Probability Measure} \label{sec:recursion}
First, we seek to characterize the backward recursion that enables propagation of probability measures that satisfy staying within a target tube for all probability values.
Satisfaction of the target tube requires additional constraints to be placed on this propagation.
To ensure that the state probability measure $\mu_{k}$ resides within the target tube set for a time step $k$, we constrain

\begin{equation}
    \mu_{k}(T_{k})= \int_{\mathbb{R}^n}\ind_{T_k}( x_k)\mu_{k}(\d x_k),
\label{eq:TargetTubeProb}
\end{equation} 
to equal $\alpha$.
For the final state probability measure to meet the target tube requirement for any probability values $\alpha\in[0,1]$, we constrain $\mu_{N-1}$,
\begin{align}\label{eq:finalProb}
    \mu_{N}(T_{N}) = \int_{\mathcal{X}}\ind_{T_N}(x_N) Q(\d x_{N}|x_{N-1})\mu_{N-1}(\d x_{N-1}),
\end{align}
to equal $\alpha$.

We utilize these constraints in the following theorem to establish satisfaction of the reach measure via state probability measures.
\begin{thm}
\label{thm:MeasureBackwardRecursion}
    If we constrain the state measures such that
    \begin{itemize}
        \item (\ref{eq:TargetTubeProb}) holds for $k \in \mathbb{N}_{[0,N-2]}$;
        \item (\ref{eq:finalProb}) holds;
    \end{itemize}
    then the system in Definition~\ref{def:TransitionKernel} satisfies the reachability specification in Problem~\ref{prob:RAviaPM} for all probability values, $\alpha$.
\end{thm}
\begin{proof}
Consider the case in which $k = N-1$:
Note that \eqref{eq:finalProb} holds  
as long as $\mu_{N-1}$ satisfies \eqref{eq:TargetTubeProb}. 
Now consider $k \in \mathbb{N}_{[0,N-2]}$.  We employ backwards induction, beginning with the base case $k = N-2$, 
\begin{align}
        &\mu_{N-1}(X)=\int_{\mathcal{X}}Q(X \mid x_{N-2})\mu_{N-2}(\d x_{N-2}),
    \end{align}
    holds $\forall X\in\mathscr{B}(\mathcal{X})$ provided that $\mu_{N-2}$ satisfies \eqref{eq:TargetTubeProb}. 
    If the equality holds for the case $k = j$ where $j < N-2$, then it must also hold for the case $k = j-1$.
    Observe that for $k = j$ and $k = j-1$ respectively, the recursions are,
    \begin{subequations}
        \begin{align}
            \mu_{j+1}(X) &=\int_{\mathcal{X}}Q(X\mid x_{j})\mu_{j}(\d x_{j}),\\
            \mu_{j}(X) &= \int_{\mathcal{X}}Q(X\mid x_{j-1})\mu_{j-1}(\d x_{j-1}),
        \end{align}
        holds $\forall X\in\mathscr{B}(\mathcal{X})$ as long as $\mu_{j}$ and $\mu_{j-1}$ satisfy \eqref{eq:TargetTubeProb}. 
        Thus, the propagation of the state probability measures satisfy the reach specification in Problem~\ref{prob:RAviaPM}.  
    \end{subequations}
\end{proof}
Theorem~\ref{thm:MeasureBackwardRecursion} draws inspiration from ~\cite[Proposition 7.28]{bertsekas1996stochastic} which proves the existence of a unique probability measure, stepping backwards in time, to an initial probability measure via the transition kernel. 
Here, we additionally constrain the state measure to reside in the target tube for some probability $\alpha$.
Nonetheless, by restricting the propagation of the state measures, we are able to characterize state probability measures addressing Problem~\ref{prob:RAviaPM}.

\subsection{Constructing Probability Measures via Level Sets}
\label{subsec:ProbMeasuresLevelSets}

The backward recursion in Theorem~\ref{thm:MeasureBackwardRecursion}, if computed directly, would yield the probability measure $\mu_{k}$ associated with the state at each time step. 
However, a direct computation constitutes a search through a space of probability measures, which is cannot practically be done in the case of continuous domains.
We instead use a more tractable formulation based on the relationship between the level sets of the probability of the state lying within the reach tube, $\alpha$, and the state probability measure.
This relationship is predicated on the relationship between the value function in standard stochastic reachability and the relationship to the measure theoretic formulation. 

To make this link clear, we reintroduce the backward recursion from~\cite{abate2008probabilistic,summers2010verification} with our reachability specification with target tubes, $T_k$, 
\begin{subequations}
    \label{eq:StandardBackwardRecursion}
    \begin{align}
V_{N-1}(x_{N-1}) &= \ind_{T_{N-1}}(x_{N-1})\mathbb{E}[\ind_{T_{N}}(\mathbf{x}_{N})| x_{N-1}],\\
V_{N-2}(x_{N-2}) &= \ind_{T_{N-2}}(x_{N-2})\mathbb{E}[V_{N-1}(\mathbf{x}_{N-1})| x_{N-2}],\\
&\hspace{0.5em}\vdots\nonumber\\
V_{k}(x_{k}) &= \ind_{T_{k}}(x_{k})\mathbb{E}[V_{k+1}(\mathbf{x}_{k+1})| x_{k}].
    \end{align}
\end{subequations}%
We compute the conditional expectation of the value function, $V_{k+1}: \mathbb{R}^n\rightarrow [0,1]$ at time step $k\in\mathbb{N}$, given $x_k$, from the transition kernel in~\eqref{eq:transitionKernel} as
\begin{equation}
    \label{eq:conditionalExpectation}
    \mathbb{E}[V_{k+1}(\mathbf{x}_{k+1})| x_k] = \int_{\mathcal{X}}V_{k+1}(x_{k+1})Q(\d x_{k+1}|x_k),
\end{equation} 
where $\mathbf{x}_{k+1}$ denotes the random state at the next time step.
To infer the probability measure from \eqref{eq:StandardBackwardRecursion} and relate it to the measure-theoretic approach, we introduce the following proposition. 
\begin{prop}
    \label{prop:equivMeasure}
    If the expectation of the value function, represented by the Lebesgue integral with respect to state probability measure $\hat \mu_k$, equals $\alpha$,
    That is
    \begin{subequations}
            \label{eq:equivStandardReach}
            \begin{align}
            &\mathbb{E}[V_{N-1}(\mathbf{x}_{N-1})] =\nonumber\\ 
            &\hspace{0.3in}\int_{\mathcal{X}}\ind_{T_{N-1}}(x_{N-1})\mathbb{E}[\ind_{T_N}(\mathbf{x}_{N})| x_{N-1}] \hat\mu_{N-1}(\d x_{N-1}), \label{eq:equivStandardReachFinal}\\
            &\hspace{0.35in}\mathbb{E}[V_{k}(\mathbf{x}_k)] = 
            \int_{\mathcal{X}}\ind_{T_{k}}(x_{k})\mathbb{E}[V_{k+1}(\mathbf{x}_{k+1})| x_{k}] \hat\mu_{k}(\d x_{k}),\nonumber\\
            &\hspace{1.9in}\mathrm{for}\ k\in\mathbb{N}_{[0,N-2]},\label{eq:equivStandardReachInter}
        \end{align}
    \end{subequations}
    both equal $\alpha$,
   then the state probability measure $\hat\mu_k$ from \eqref{eq:equivStandardReach} and $\mu_k$ from Theorem~\ref{thm:MeasureBackwardRecursion} are equivalent.
\end{prop}
\begin{proof}
    We show this by stepping backwards in time, starting with $k=N-1$.
    Given $\hat\mu_{N-1}$ which satisfies \eqref{eq:equivStandardReachFinal}, it also satisfies both \eqref{eq:TargetTubeProb} and \eqref{eq:finalProb} for some $\alpha$.
    Therefore, $\hat\mu_{N-1} = \mu_{N-1}$.
    Likewise, given state probability measures $\hat\mu_{k}$ which satisfy \eqref{eq:equivStandardReachInter} for $k\in\mathbb{N}_{[0,N-2]}$, they must satisfy \eqref{eq:TargetTubeProb} since, 
    \begin{equation}
        \hat\mu_{k+1}(X) = \int_{\mathcal{X}}Q(X\mid x_k)\hat\mu_k(\d x_k),
    \end{equation}
    $\forall X\in\mathscr{B}(\mathcal{X}) $via \eqref{eq:measureState}, results in
    \begin{equation}
        \hat\mu_{k+1}(T_{k+1}) = \int_{\mathcal{X}}\ind_{T_{k+1}}(x_{k+1})\hat\mu_{k+1}(\d x_{k+1}), 
    \end{equation}
    equaling $\alpha$.
    Thus $\hat\mu_{k} = \mu_{k}$.
\end{proof}
The value function determines if the current state, $x_k$, is within the target tube prior to evaluating the conditional expectation, \eqref{eq:conditionalExpectation}. 
This conditional expectation represents both the state probability measure via the transition kernel and whether each state probability measure lies in the reach tube for the rest of the time horizon.
In contrast, as shown in Theorem~\ref{thm:MeasureBackwardRecursion}, our approach constrains the current state probability measure, $\mu_k$, \textit{prior} to the propagation to the next constrained state probability measure, $\mu_{k+1}$, which we define via the transition kernel.
More importantly, what Proposition~\ref{prop:equivMeasure} identifies is that, for each time step \(k\in\mathbb{N}_{[0,N-1]}\), the evaluation of the value function at state $x_k$ corresponds to the level sets of the state probability measure. 
That is, by identifying the states $x_k$ that align with a given level set probability, $\alpha$, we can infer the respective level sets of the state probability measures, $\mu_k$, such that \eqref{eq:equivStandardReach} holds.

\section{Approximating the Level Sets of the State Probability Measure}
\label{sec:NN}

To obtain the exact probability measures $\mu_k$ would require a direct solution of~\eqref{eq:StandardBackwardRecursion} using direct knowledge of the transition kernel. In cases where this knowledge is not available, we must instead turn to an approximation using what information is available. In this section, we consider an approximation method based on the assumption that the only information available to us is trajectory data. 
Our ultimate goal is to form a surrogate of the measure as shown in Figure~\ref{fig:functionapproximator}: a function approximation that returns, for a given state value and time step, as accurate an estimate of $V_k$ as our data allow.

The essence of our strategy is as follows: first, form empirical estimates $\hat{V}_k$ of the measures $\mu_k$ using sample data; then, use the empirical estimates to form training data comprising (state, time step, empirical measure) tuples; and finally, to train a function approximator on that data.

\subsection{Forming the Training Data}

We assume that the information available to us takes the following form.

\begin{assum}
    \label{assum:samplingState}
    The transition kernel is not given, but we can compute state samples, i.e. $x_{k+1,i}\sim Q(\cdot|x_k)$ where $x_{k+1,i}\in\mathbb{R}^n$ provided some $x_k\in\mathbb{R}^n$,
    where $x_{k}\in\mathbb{R}^n$ is the current state and $x_{k+1,i}\in\mathbb{R}^n$ is the sample of the next state.
\end{assum}

From these data we compute the empirical estimates
\begin{subequations}
    \label{eq:StandardBackwardRecursionEmpirical}
    \begin{align}
    \hat{V}_{N-1}(x_{N-1}) &= \ind_{T_{N-1}}(x_{N-1})\frac{1}{L}\sum_{j=1}^L\ind_{T_{N}}(x_{N,j}),\\
    \hat{V}_{N-2}(x_{N-2}) &= \ind_{T_{N-2}}(x_{N-2})\frac{1}{L}\sum_{j=1}^L \hat{V}_{N-1}(x_{N-1,j}),\\
    &\hspace{0.5em}\vdots\nonumber\\
    \hat{V}_k(x_{k}) &= \ind_{T_{k}}(x_{k})\frac{1}{L}\sum_{j=1}^L \hat{V}_{k+1}(x_{k+1,j}).
    \end{align} 
\end{subequations}
Here, $L$ denotes the number of samples of the subsequent states $x_{k+1,j}$ starting from a state $x_k$ that is sampled from a uniform distribution over the state space, i.e. $x_k\sim\mathrm{Uniform}(\mathcal{X})$. 
We restrict our attention to a hyperrectangle in the state space, taking $\mathcal{X} = [-b,b]^n$ with bounds $b\in\mathbb{R}$, which we assume encompass the tube sets $T_k$.
This process is contingent on the fact that since we are handling dynamical systems which have no control input, we can avoid a backward recursion via dynamic programming and, instead, employ forward simulation to obtain an empirical estimate.
We also note that this process would not be effective for a dynamical system with a control input.

The following algorithm outlines the procedure to generate the sample trajectories $\{\{x_{l,j}\}_{j=1}^L\}_{l=k+1}^{N}$ given a single state sample $x_k$ and $L$ $x_{k+1}$ samples via Assumption~\ref{assum:samplingState}. 
\begin{algorithm}[ht!]
\caption{Trajectory generation procedure.}\label{alg:sampling}
    \begin{algorithmic}[1]
        \Require $x_k\in\mathbb{R}^n$, $L\in\mathbb{N}$
        \Ensure $\{\{x_{l,j}\}_{j=1}^L\}_{l=k+1}^{N}$ 
        \For{$l \in \mathbb{N}_{[k,N-1]}$}
        \For{$j\in \mathbb{N}_{[1,L]}$}
        \State $x_{j,l+1}\sim Q(\cdot|x_k)$ 
        \EndFor
        \EndFor
    \end{algorithmic}
\end{algorithm}

To demonstrate the consistency of our choice of empirical measures, the following result shows that \eqref{eq:StandardBackwardRecursionEmpirical} converges to \eqref{eq:StandardBackwardRecursion} in the limit for a single evaluation point of the state, $x_k$, as $L$ goes to infinity. 
\begin{prop}
    \label{prop:convergenceValue}
    If the state samples are produced according to Algorithm~\ref{alg:sampling}, then $\lim_{L\rightarrow\infty} \hat{V}_k(x_{k}) = V_k(x_k),\ \forall k\in\mathbb{N}_{[0,N-1]}$.   
\end{prop}
\begin{proof}
    Let~\eqref{eq:StandardBackwardRecursion} denote the true value function that is defined via a Lebesgue integral~\cite[Ch. 2, Definition 2.4.1]{bogachev2007measure} of a function with respect to the previous transition kernels as well as prior state probability measures and let~\eqref{eq:StandardBackwardRecursionEmpirical} denote an empirical value function averaged over $L$ samples detailed in Algorithm~\ref{alg:sampling}.
    The strong law of large numbers~\cite[Ch. 10, 10.10(v)]{bogachev2007measure} ensures, as $L$ increases, $\hat{V}_k$ will converge almost everywhere to $V_k$.
    Thus, for $k\in\mathbb{N}_{[0,N-1]}$,
    $\lim_{L\rightarrow\infty}\hat{V}_k(x_k) = V_k(x_k)$.
\end{proof}

Having established the trajectory data and empirical measures, we proceed to the second step of our strategy: forming the training data.
We associate to each trajectory datum $x_k$ a value of the corresponding measure value as estimated by the empirical measures $\hat{V}_k$ to form these estimates. Algorithm~\ref{alg:sampling} demonstrates the exact procedure.
\begin{algorithm}[h!]
\caption{Training data generation procedure.}\label{alg:trainData}
\begin{algorithmic}[1]
\Require $[-b,b]^n\subseteq\mathbb{R}^n$, $M,L\in\mathbb{N}$
\Ensure $\mathrm{data}\in\mathbb{R}^{M\times N}$
\For{$k \in \{N-1, N-2,\hdots,1,\ 0\}$}
\For{$i\in\mathbb{N}_{[1,M]}$}
    \State $x_i\sim\mathrm{Uniform}([-b,b]^n)$
    \State $\{\{x_{l,j}\}_{j=1}^L\}_{l=k+1}^{N}$ via Algorithm~\ref{alg:sampling} 
    \State $\hat{\alpha}_{i,k} = \hat{V}_k(x_{k,i})$ with $\{\{x_{l,j}\}_{j=1}^L\}_{l=k+1}^{N}$.
    \State $\mathrm{data}_{i,k} = ((x_{i,k},k),\hat{\alpha}_{i,k})$
\EndFor
\EndFor
\end{algorithmic}
\end{algorithm}
\begin{figure}
    \centering
    \tikzset{
    every picture/.style={line width=0.75pt}, 
    arrow/.style={-{Triangle[length=3mm, width=2mm]}, color=black, line width=0.75pt}
}

\begin{tikzpicture}[x=0.75pt,y=0.75pt,yscale=-1,xscale=1]
    \draw [arrow] (250.33,159.77) -- (290.33,159.77);
    \draw [arrow] (250.33,189.77) -- (290.33,189.77);
    
    \draw   (291,148.77) rectangle (431.33,199.77);
    
    \draw [arrow] (431.33,175.77) -- (468.33,175.77);
    
    \draw (361.67,174.27) node [anchor=center] {Function Approximator};
    \draw (250,159.17) node [anchor=east] {$x$};
    \draw (250,189.77) node [anchor=east] {$k$};
    \draw (470,175.77) node [anchor=west] {$\alpha$};
\end{tikzpicture}
    \caption{The function approximator takes as input the state and time step, then outputs a reach probability value of satisfying the reach specification.}
    \label{fig:functionapproximator}
\end{figure}
To facilitate the training process, we restructure the training dataset, denoted by $\mathrm{data}$, into $\mathrm{data}'$ so that each entry, $\mathrm{data}'_{i'}$, corresponds to $\mathrm{data}_{i,k}$ for every $i$ within $[1,M]$ and $k$ within $[0,N-1]$. 
Thus, the training dataset comprises a collection of tuples, $\mathrm{data}' = \{(x_{i'},k_{i'}), \hat{\alpha}_{i'}\}_{i'=1}^{M+N}$, where the input includes the state and the time step, while the output is the probability value corresponding to that state and time step.

With the training data formed, we proceed to the final step: training the function approximator.
Training a function approximator, $\hat{V}_\theta: \mathcal{X}\times \mathbb{N} \rightarrow [0,1]$ according to the training data tuples from the previous step corresponds to solving the following optimization problem:
\begin{align}
    \label{opt:training}
    \underset{\theta}{\mathrm{minimize}} &\quad \frac{1}{M+N}\sum_{i'=1}^{M+N}J\left(\hat{V}_{\theta}(x_{i'},k_{i'}),\hat{\alpha}_{i'}\right),
\end{align}
where $J: \mathbb{R}\times \mathbb{R}\rightarrow\mathbb{R}$ is a loss function and $\theta$ represents the parameters of the function approximator.
Once the training is complete, we may use the resulting approximator to solve Problem~\ref{prob:RAviaPM} for any reachability specification we like: the procedure reduces to numerical root-finding to learn the level sets of the state probability measure that satisfy the reach tube specification at time step $k$.

\begin{figure*}[t!]
  \centering
  \subfloat[Double integrator with variance of $\sigma^2 = 0.01$]{\includegraphics[width=\linewidth]{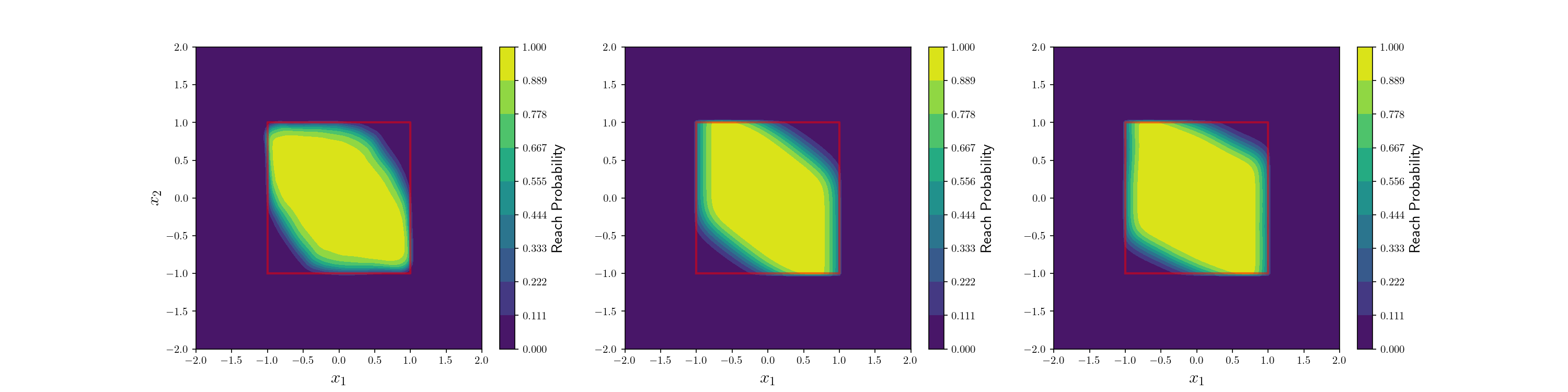}
  \label{fig:GaussDI1}}
  
  \subfloat[Double integrator with variance of $\sigma^2 = 0.1$]{\centering\includegraphics[width=\linewidth]{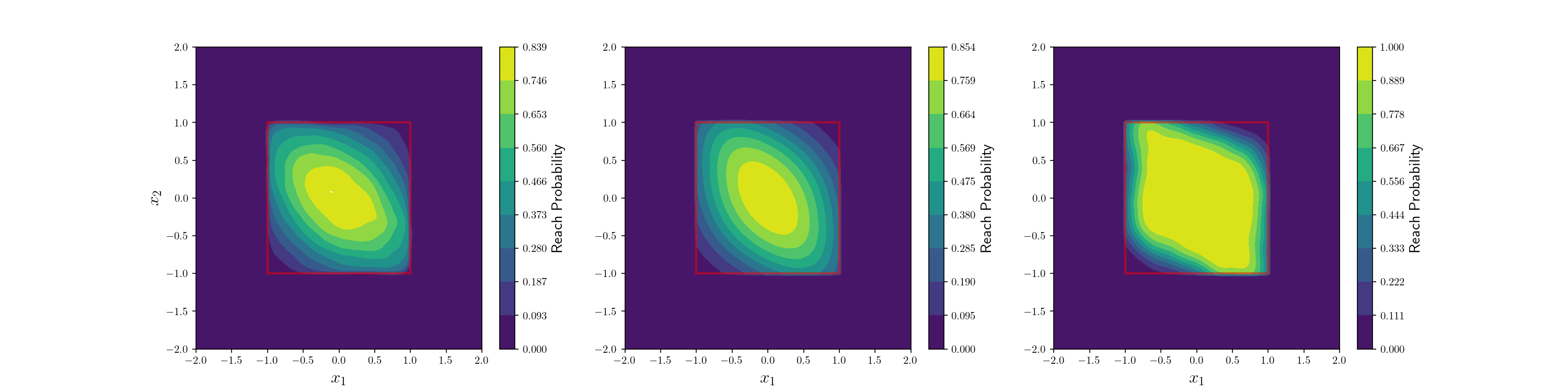}
  \label{fig:GaussDI2}}
  \caption{We compare our approach to the ground truth via dynamic programming (DP) as well as to a reproducing kernel Hilbert space (RKHS) approach~\cite{thorpe2022state}.
  In the top row, for a Gaussian disturbance with $\sigma^2 = 0.01$, the RKHS (right) outperforms our method (left), matching closely with DP (middle), as the neural network struggles to learn the sharp drop off of probabilities at the boundaries.
  However, in the second row, when the Gaussian disturbance has $\sigma^2 = 0.1$ variance, our approach (left) is close to the DP solution (middle) in comparison to RKHS (right).
  Note that we are able to train the neural network with more samples ($M = 5000$ state samples at training time with $L=2000$ to compute $\alpha$ offline) versus the RKHS method ($M = 10,201$ samples at training time). 
  When attempting to use the same number of samples ($M\cdot L = 10,000,000$ samples) for the RKHS approach, we ran out of memory.}\label{fig:SRDoubleIntegrator}
\end{figure*}

\subsection{Using DNNs as the Function Approximator}

With the strategy for training a function-approximator surrogate for the measures established for an abstract function approximator, we turn to the design choices required to implement the method using a DNN. 
There are three design choices to make: the loss function, the layer activation functions, and the layer architecture. 

\textbf{Activation Functions}: Since the output must lie between zero and one, the last layer uses a sigmoid function. All others use ReLU functions.

\textbf{Loss Function}:
Recall that the training problem at hand is to estimate a measure from point evaluations. For this type of supervised learning problem, an effective loss function is the binary cross entropy (BCE) loss, 
\begin{align}
    BCE &= -\cfrac{1}{s}\sum_{i'=1}^s\Bigg[\hat{\alpha}_{i'} \cdot \log\left(\Tilde{V}_{\theta}(x_{i'},k_{i'})\right) \nonumber \\
    &\hspace{4em} +  (1-\hat{\alpha}_{i'})\cdot \log\left(1-\Tilde{V}_{\theta}(x_{i'},k_{i'})\right)\Bigg].
\end{align}
Here, $s\in\mathbb{N}$ denotes the batch size of the data that is being learned in a single pass of the training loop, which continues until the model has trained on the entire data set, i.e. $s<M+N$, then repeats the process over again for a finite number of iterations known as epochs.

\textbf{Neural Network Structure}:
We design our DNN as a fully-connected, feed-forward neural network. 
We use $4$ hidden layers where each layer consists of $64$ neurons, each employing the Rectified Linear Unit (ReLU) activation function. 
The input into the neural network is the state, $x\in\mathcal{X}$, and the current time step, $k\in\mathbb{N}$.

Before turning to the examples, we briefly consider the matter of computational complexity involved in evaluating the function approximator.
Note that DNNs evaluations scale with the architecture's size, specifically the arrangement and connectivity of its layers, and not the sample size.
This means that a trained neural network can accommodate additional data without a proportional increase in the evaluation time.
In comparison, kernel methods, which scale with the number of samples, struggle at evaluation time due to the size of the Gram matrix~\cite{drineas2005nystrom}.
Although kernel methods have seen scalability improvements via random Fourier features~\cite{rahimi2007random}, neural networks have made parallel strides in scalability through techniques such as quantization~\cite{nagel2021white}, to represent neural networks on a smaller memory footprint, and batching of the data for training with limited memory~\cite[Ch. 8]{Goodfellow-et-al-2016}.
    
\section{Examples}
\label{sec:example}
We demonstrate our approach on three examples. 
For training the algorithm we use PyTorch~\cite{paszke2019pytorch} with the Adam optimizer~\cite{kingma2014adam}. 
Specifically, we use a learning rate scheduler, which initially starts with a rate of $0.001$, and adjusts the learning rate by multiplying it by $0.1$ when the loss function flattens, i.e. stops improving, during the training.
We use an Xavier initialization scheme to prevent exploding and vanishing gradients which keeps the variance of the activation functions the same across each layer~\cite{glorot2010understanding}.
All computations were done in Python on an Apple M1 Macbook Pro with 16GB of RAM with PyTorch running in CPU mode, making no use of GPU or neural network hardware acceleration. 
For comparisons to the proposed approach, we employ both dynamic programming (DP)~\cite{SReachTools} and a reproducing kernel Hilbert spaces (RKHS) approach~\cite{thorpe2022state}.

\subsection{Double Integrator Experiments}
\subsubsection{2D Double Integrator}
This example compares our approach against dynamic programming and RKHS \cite{thorpe2022state} approaches.
Consider dynamics of a double integrator with time horizon $N=3$ and sampling time $\Delta T = 0.25$,
\begin{equation}
    \mathbf{x}_{k+1} = \begin{bmatrix} 1 & \Delta T \\ 0 & 1 \end{bmatrix}\mathbf{x}_k + \mathbf{w}_k,
\end{equation}
where $\mathbf{x}_{k}\in\mathcal{X} = [-2,2]^2$ is the random state vector. 
We consider Gaussian disturbances, $\mathbf{w}_{k}\sim\mathcal{N}(0_2,\sigma^2 I_{2\times 2})$.
We consider variances of $\sigma^2 = 0.01$ and $\sigma^2 = 0.1$.
The target sets are $T_k=[-1,1]^2$ for $k\in\mathbb{N}_{[0,N]}$.

We trained the DNN with $M=5000$ initial samples of the state and $L=2000$ state samples, which we split into batches of size $101$ over $50$ epochs to train the neural network. 
Results for both values of $\sigma$ are shown in Figure \ref{fig:GaussDI1}, as well as comparisons with RKHS and dynamic programming. 
We train the RKHS approach with a sample size of $10,201$ for both variances. 
For the example with variance $\sigma^2=0.01$, the RKHS approach 
matches closely with the dynamic programming solution, while our approach has error 
at drop offs in probability at the edges (Figure~\ref{fig:GaussDI1}).
However, with variance $\sigma^2=0.1$, the RKHS approach does not accurately match the dynamic programming solution.
In contrast, our approach is able to handle the larger noise, and is close to the dynamic programming solution, as seen in Figure~\ref{fig:GaussDI2}. 
We report the errors for both examples in Table~\ref{table:computetimeerror}.
When using the same number of state samples ($M\cdot L = 10,000,000 $) as we use to train the DNN, we ran out of memory attempting to train the RKHS method.

\subsubsection{$n$-Dimensional Stochastic Chain of Integrators}
\begin{figure}[t]
    \centering    \includegraphics[height=2.05in,width=\linewidth]{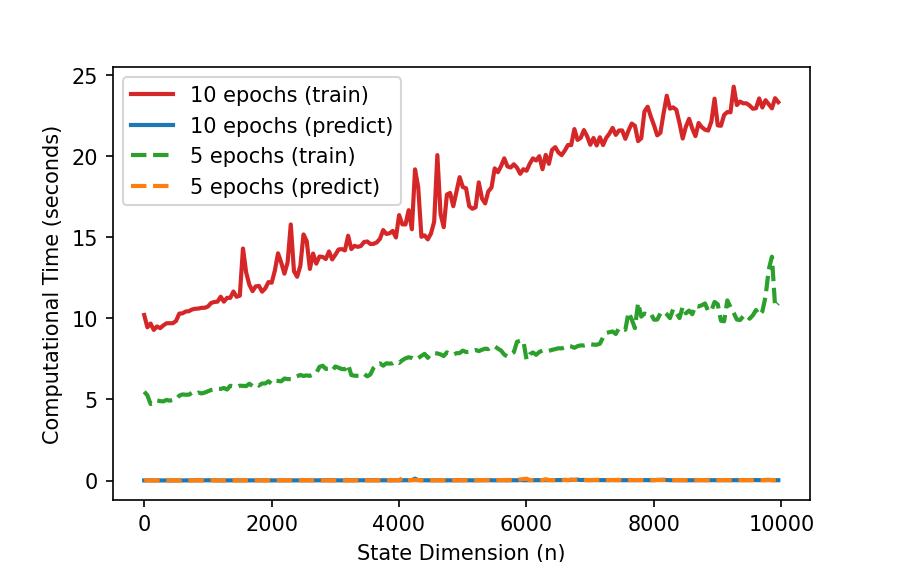}
    \caption{This plot shows that the training time scales with dimension for the neural network when we fix the number of samples during training.
    For a fixed sample size, but with increasing state dimension, the neural network scales linearly, similarly to RKHS~\cite{thorpe2019model}.
    Note that increasing the number of epochs increases the slope of the training curve.}
    \label{fig:nDim}
\end{figure}

In this example, we show that our approach scales linearly with dimension.
We consider the $n$-dimensional stochastic chain of integrators from~\cite{vinod2017scalable,thorpe2019model}, 

\begin{align}
    \mathbf{x}_{k+1} = 
    \begin{bmatrix}
        1 & \Delta T & \cdots & \frac{\Delta T^{n-1}}{(n-1)!} \\
        0 & 1 & \cdots & \frac{\Delta T^{n-2}}{(n-2)!} \\
        \vdots & \vdots & \ddots & \vdots \\
        0 & 0 & \cdots & 1
        \end{bmatrix} \mathbf{x}_k  + \mathbf{w}_k,
\end{align}
where we vary $n$ from $2$ to $10,000$ dimensions, $\mathbf{x}_{k}\in\mathcal{X}\subseteq\mathbb{R}^n$ is the random state vector evaluate with a zero vector $0_{n}$, and $\mathbf{w}_{k}\sim\mathcal{N}(0_n,0.01I_{n\times n})$ is the Gaussian noise vector, sampling rate is $\Delta T = 0.25$, and the time horizon is $N=5$.
For the reachability specification, we specify the tube sets to be $T_k=[-1,1]^n$ for $k\in\mathbb{N}_{[0,N]}$.
We trained our DNN with $M=1024$ initial samples of the state and $L=1$ state samples, which we split into batches of $10$ samples during training time. 
We iterate through the entire training set for $5$ and $10$ epochs, respectively, to determine how the training procedure scales with additional passes through the training set. 
Our approach scales linearly and has performance similar to RKHS~\cite{thorpe2019model}, as shown in Figure \ref{fig:nDim}.
In contrast to RKHS, in which samples are contained within the Gram matrix, a neural network can utilize batch processes to train over large amounts of data, without increases to the time needed to evaluate the neural network.
Note that the increase in epochs results in a steeper slope.

\subsection{Quaternion Attitude Dynamics}
This example demonstrates solution to stochastic reachability problems with nonlinear dynamics.
The rigid body dynamics of a torque-free, tumbling object is given by Euler's rotation equations,
\begin{align} \label{eq:eulerrotation}
        \dot{\omega}_x &= \left[(I_{yy} - I_{zz})/ I_{xx}\right] \omega_y \omega_z\nonumber \\
        \dot{\omega}_y &= \left[(I_{zz} - I_{xx})/I_{yy}\right] \omega_z \omega_x, \\
        \dot{\omega_z} &= \left[(I_{xx}- I_{yy})/I_{zz}\right] \omega_x \omega_y\nonumber
\end{align}
where $I\in\mathbb{R}^3$ is the inertia matrix with diagonal entries, $I_{xx}$, $I_{yy}$, and $I_{zz}$, since we assume no products of inertia, and angular velocity about the principle axes is denoted by $\omega = [\omega_x \ \omega_y \ \omega_z]^\top \in \mathbb{R}^{3}$.
The evolution of a unit quaternion representing attitude in the earth-centered inertial frame with the angular velocity in the body frame is
\begin{align} \label{eq:quat}
    \dot{q} =
    \frac{1}{2}
    q
    \otimes 
    \begin{bmatrix}
    0 \\
    \mathbf{\omega} \\
    \end{bmatrix},
\end{align}
where $q = a + bi + cj + dk$ is a unit quaternion on the $3$-sphere, denoted by $\mathbb{S}^3$, and it is of unit length, i.e., $a^2+b^2+c^2+d^2=1$~\cite{jia2008quaternions}.
The state vector is $x = [\omega_x,\omega_y,\omega_z,q_w,q_x,q_y,q_z]^\top\in\mathbb{R}^7$, where $[\omega_x, \omega_y,\omega_z]^\top$ represents the angular velocity about the principal axes, $q_w$ is a scalar that represents the angle of rotation, and $[q_x,q_y,q_z]^\top$ is a unit vector that represents the axis of rotation.

To obtain~\eqref{eq:transitionKernel}, we discretized~\eqref{eq:eulerrotation} and~\eqref{eq:quat} via a fourth-order Runge-Kutta scheme, then 
added Gaussian noise with mean $0$ and variance $\sigma^2=0.1$ to each of the elements of the quaternion kinematics.
We generated points on the surface of the unit sphere, which we represent as pure quaternions, i.e. $q=[0,x,y,z]$, then simulated the attitude dynamics via fourth-order Runge Kutta with time horizon $N=3$, sampling time $\Delta T = 0.25$, initial constant angular velocity $\omega_0 = [0,0.5,1]^\top$, and moments of inertia $I_{xx} = 10$, $I_{yy}=5$, and $I_{zz}=7$. 
Our reachability specification is on a projected space comprised of elevation and azimuth, i.e. $T_k=\left\{\theta, \phi \mid \frac{\pi}{4} \leq \theta \leq \frac{\pi}{2}, \frac{\pi}{4} \leq \phi \leq \frac{3\pi}{4}\right\}$ for $k\in\mathbb{N}_{[0,N]}$, which we can derive via the unit vector component of the quaternion that represents a point in three-dimensional space.
\begin{subequations}
    \begin{align}
        \theta &= \arccos(q_z)\\
        \phi &= \mathrm{arctan}(q_y/q_x)
    \end{align}
\end{subequations}%
These equations are standard for converting unit quaternions to spherical coordinates~\cite{marsden2003vector}.

\begin{figure}[ht!]
  \centering
  \subfloat[Our approach]{\includegraphics[width=\linewidth]{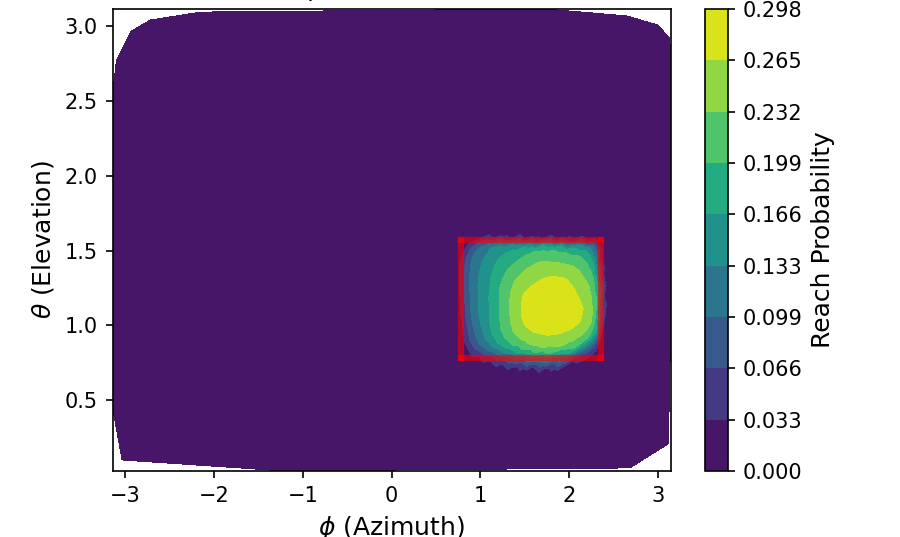}
  \label{fig:DNNQuaternion}}
  \vspace{0.25cm}
  \subfloat[Empirical estimate via~\eqref{eq:StandardBackwardRecursionEmpirical}]{\includegraphics[width=\linewidth]{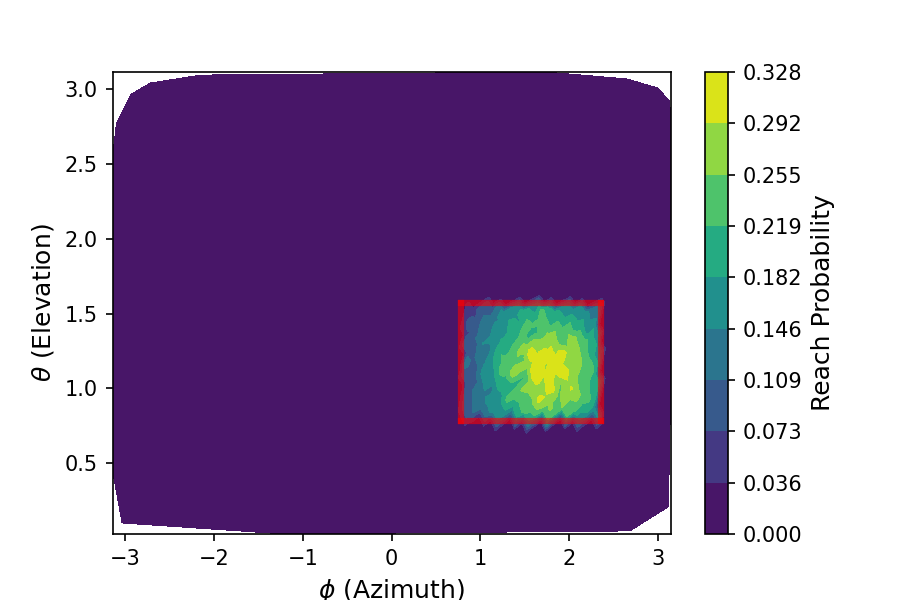}
  \label{fig:EmpiricalQuaternion}}
  \vspace{0.25cm}
  \subfloat[RKHS approach]{\includegraphics[width=\linewidth]{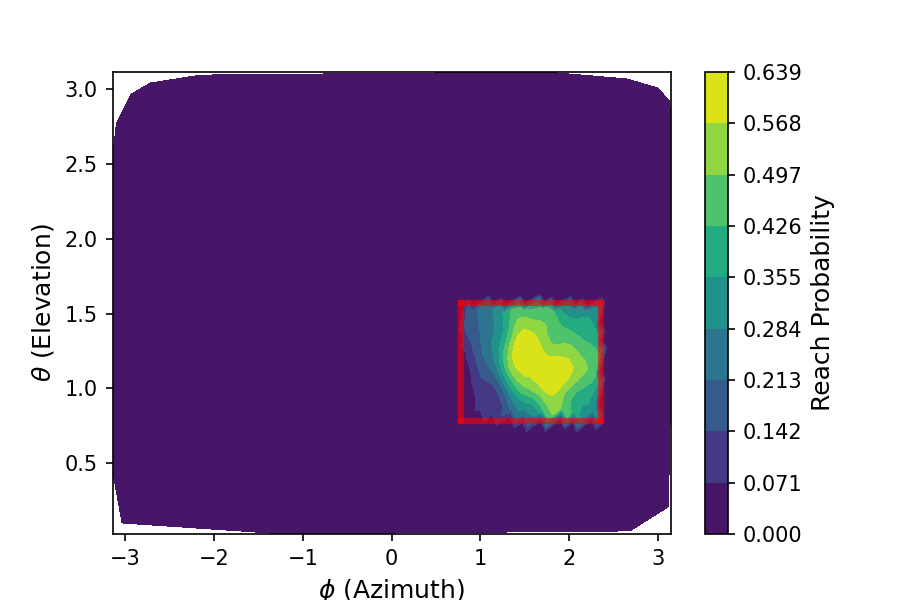}
  \label{fig:RKHSQuaternion}} 
  \caption{We compare the proposed approach with an empirical estimate via~\eqref{eq:StandardBackwardRecursionEmpirical}, and the RKHS approach.
  The neural network (Figure \ref{fig:DNNQuaternion}) is closer to the empirical estimate (Figure \ref{fig:EmpiricalQuaternion}) than the RKHS (Figure~\ref{fig:RKHSQuaternion}).
  }
  \label{fig:attitude}
\end{figure}

We trained our DNN with $M=2500$ initial samples of the state and $L=500$ state samples, which we split into batches of size $10$ over $50$ epochs to train the neural network. 
Due to the dimensionality of the state, we could not utilize dynamic programming. 
However, as shown in Figure \ref{fig:attitude}, the DNN approach matches closely with the empirical estimate via~\eqref{eq:StandardBackwardRecursionEmpirical}, whereas the RKHS does not match closely and yields higher probability values. 
The error with respect to~\eqref{eq:StandardBackwardRecursionEmpirical} is provided in Table ~\ref{table:computetimeerror}.  The RKHS approach has both a higher maximum absolute error and average error, 
which could indicate that it requires additional hyperparameter tuning, cross-validation, and additional data.

The differences in runtime between the RKHS and DNN approaches in the double integrator and attitude dynamics examples are influenced by the complexity of the systems and number of samples used. In the double integrator example, the RKHS takes longer to compute due to the cost of evaluating a large Gram matrix with more samples. However, in the attitude dynamics example, the RKHS method is quicker, likely because the kernel evaluations and matrix operations are less computationally intensive for this particular sample size of 2500 compared to the DNN’s extensive training requirements.

\begin{table}[ht]
    \centering
    \textbf{Double Integrator with variance} $\boldsymbol{\sigma^2=0.01}$
    \begin{tabular}{lccccc}
    \hline
         & State & Max Absolute & Average & \\
         Method & Samples & Error with DP & Error with DP & Time (s)  \\
    \hline
         DNN & $5000$ & $0.8433$ & $0.0205$ & $18.905$ \\
         DP & $10201$ & -- & -- & 8.99 \\
         RKHS & $10201$ & $0.4074$ & $0.0142$ & $53.23$ \\
    \hline\\
    \end{tabular}
    
    \textbf{Double Integrator with variance} $\boldsymbol{\sigma^2=0.1}$
    \begin{tabular}{lccccc}
    \hline
         & State & Max Absolute & Average & \\
         Method & Samples & Error with DP & Error with DP & Time (s)  \\
    \hline
         DNN & $5000$ & $0.3106$ & $0.0063$ & $20.35$ \\
         DP & $10201$ & -- & -- & $9.13$ \\
         RKHS & $10201$ & $0.5091$ & $0.0669$ & $50.158$ \\
    \hline\\
    \end{tabular}
    
    \textbf{Quaternion Attitude Dynamics}
    \begin{tabular}{lccccc}
    \hline
         & State & Max Absolute & Average & \\
         Method & Samples & Error with \eqref{eq:StandardBackwardRecursionEmpirical} & Error with \eqref{eq:StandardBackwardRecursionEmpirical} & Time (s)  \\
    \hline
         DNN & $2500$ & $0.1477$ & $0.0023$ & $52.186$ \\
         RKHS & $2500$ & $0.3901$ & $0.0164$ & $1.108$ \\
    \hline\\
    \end{tabular}
    \caption{Number of state samples, maximum absolute error, average error, and computation times for each method for the double integrator and the quaternion attitude dynamics.
    The time reported for dynamic programming consists of evaluation on a uniform grid to compute the probabilities via backward recursion.
    For the DNN and RKHS approaches, we show the combined training and evaluation times.}
    \label{table:computetimeerror}
\end{table}

\section{Conclusion}
\label{sec:conc}

We provided a measure theoretic formalism of the stochastic reachability problem for uncontrolled systems, that characterizes stochastic reachability in terms of probability measures.  This approach enables new computational tools for approximate solutions to stochastic reachable sets, by computing level sets of the probability measure for all probability values. 
We developed a numerical implementation that employs deep neural nets to approximate the stochastic reachability probability, for a given state and time step within a finite time horizon.  
Future work will focus on extensions to accommodate dynamical systems with control inputs. The essence of this extension will be to generalized the transition kernel framing of the dynamics into one that admits inputs, and to extend the sampling procedure to sample over inputs as well. This will require a restriction of the space of control inputs, at least initially, to a finite dimensional space, e.g. to piecewise constant signals.

\section*{Acknowledgements}
We would like to thank Adam J. Thorpe and Kendric R. Ortiz for providing their code to replicate the results for the double integrator example in~\cite{thorpe2022state} and for several discussions.
We also thank Krishna C. Kalagarla for many discussions regarding our numerical implementation via deep neural networks.

\bibliographystyle{IEEEtran}
\bibliography{refs}

\end{document}